\documentclass[12pt]{amsart}
\usepackage{amsmath,amssymb,amsthm,amscd}
\makeatletter \@addtoreset{equation}{section}

\makeatother \makeatletter

\def \<{\langle}
\def \>{\rangle}

\def \a{\alpha }
\def \e{\epsilon }
\def \l{\lambda }

\def \b{\beta }

\newtheorem{theorem}{Theorem}[section]
\newtheorem{corollary}{Corrolary}[section]
\newtheorem{lemma}{Lemma}[section]

\newtheorem{remark}{Remark}[section]
\newtheorem{definition}{Definition}[section]
\newtheorem{proposition}{Proposition}[section]

\newcommand{\bea}{\begin{eqnarray}}
\newcommand{\eea}{\end{eqnarray}}
\newcommand{\be}{\begin {equation}}
\newcommand{\ee}{\end{equation}}
\newcommand{\g}{\frak g}
\newcommand{\hg}{\hat {\frak g} }

\newcommand{\h}{\frak h}

\newcommand{\Z}{\Bbb Z}
\newcommand{\N}{\Bbb N}
\newcommand{\C}{\Bbb C}

\newcommand{\1}{\bf 1}
\newcommand{\la}{\langle}
\newcommand{\ra}{\rangle}

\def \l {\lambda}
{\vskip-\lastskip\medskip
  \noindent
  {\em #1.}\enspace
  }%
{\qed\par\medskip
  }

\def \l {\lambda}
\def \a {\alpha}
\begin{document}
\title{ Regularity of certain vertex operator algebras with
two generators }
\author{ Dra\v zen Adamovi\' c }
\date{}

\address{ Department of Mathematics, University of Zagreb,  \\
Bijeni\v cka 30, 10 000 Zagreb, Croatia} \pagestyle{plain}

  {} \maketitle

\begin{abstract}
For every $m \in {\C} \setminus \{ 0, -2\}$ and every nonnegative
integer $k$ we define the vertex operator (super)algebra $D_{m,k}$
having two generators and rank $ \frac{3 m}{m + 2}$. If $m$ is
 a positive integer then $D_{m,k}$ can be realized  as a subalgebra
of a lattice vertex algebra. In this case, we prove that $D_{m,k}$
is a regular vertex operator (super)algebra and find the number of
inequivalent irreducible modules.
\end{abstract}

\section{Introduction}

In the theory of vertex operator (super)algebras, the
classification and the construction of rational vertex operator
(super)algebras are main problems. The rationality of certain
familiar vertex operator (super)algebras was proved in papers
\cite{W}, \cite{FZ}, \cite{D}, \cite{DL}, \cite{Li1}, \cite{A1},
\cite{A2}, \cite{A3}.

It is natural to consider rational vertex operator (super)algebras
of certain  rank. In particular, in the rank one case for every
positive integer $k$ we have the well-known rational vertex
operator (super)algebra $F_k$  associated to the  lattice
$\sqrt{k}{\Z}$. These vertex operator (super)algebras are
generated by two generators.

In the present paper we will be concentrated on vertex operator
(super)algebras of rank $c_m = \frac{ 3 m} { m + 2}$, $m \in {\C }
\setminus \{0, -2 \} $. This rank has the vertex operator algebra
$L(m,0)$ associated to the irreducible vacuum $\hat{
sl_2}$--module of level $m$ and the vertex operator superalgebra
$L_{ c_m}$ associated to the  vacuum module for the $N=2$
superconformal algebra with central charge $c_m$ ( cf. \cite{A2},
\cite{A3}, \cite{FST1}, \cite{FST2}, \cite{EG}, \cite{HM}).  In
the case $m=1$ these vertex operator (super)algebras are included
into the family $F_{k}$, $k \in {\N}$, since $L(1,0) \cong F_2$
and $L_{c_1} \cong F_3$. The main purpose of this article is to
include $L(m,0)$ and $L_{ c_m}$ into the family $D_{m,k}$, $k \in
{\N}$,  of rational vertex operator (super)algebras of rank $c_m$
for arbitrary positive integer $m$.

In fact, for every $m \in {\C} \setminus \{ 0, -2 \} $ we define
the
 vertex operator (super)algebra $D_{m, k}$ as a subalgebra of the vertex
 operator (super)algebra $L(m,0) \otimes F_k$ (cf. Section \ref{def-d}). In the special case $k=1$,
  $D_{m,1}$ is in fact the  $N=2$ vertex
operator superalgebra $L_{ c_m}$ constructed  using the
Kazama-Suzuki mapping (cf. \cite{KS}, \cite{FST1}).  We also have
that $D_{m,0} \cong L(m,0)$ and $D_{1,k} \cong F_{k+2}$.

 Moreover, we
shall  demonstrate that $D_{m,k}$ has many properties similar to
affine  and $N=2$ superconformal vertex  algebras.

When $m$ is not a nonnegative integer, then $D_{m,k}$ has
infinitely many irreducible representations. Thus, it is not
rational (cf. Section \ref{def-d}). In order to construct new
examples of rational vertex operator (super)algebras we shall
consider the case when $m$ is a positive integer. Then $D_{m,k}$
can be embedded into a lattice vertex algebra (cf. Section
\ref{construction}).  In fact, we shall prove that
\bea
&& \label{par}D_{ m , k} \otimes F_{ - k} \cong L(m, 0) \otimes
F_{ - \frac{k}{2} ( m k + 2)} , \quad ( k \ \mbox{even} ) , \\
&& \label{nepar}D_{ m , k} \otimes F_{ - k} \cong L(m, 0) \otimes
F_{ - {k} ( m k + 2)} \oplus L(m,m) \otimes MF_{ - {k} ( m k + 2)}
, \quad ( k \ \mbox{odd} ) . \eea
These relations completely determine the structure of $D_{ m , k}
\otimes F_{ - k}$ as a weak $L(m,0)$--module.

 In \cite{DLM} was
introduced the notion of  regular vertex operator algebra, i. e.
rational vertex operator algebra with the property that every weak
module is completely reducible. The relations (\ref{par}) and
(\ref{nepar}), together with the regularity results from
\cite{DLM} and \cite{Li2}  imply that $D_{m,k}$ is a regular
vertex operator algebra if $k$ is even, and a regular vertex
operator superalgebra if $k$ is odd.

Let us here discuss  the case $k = 2 n$, where $n$ is a positive
integer. The relation (\ref{par}) suggests that one can study the
dual pair $(D_{m, 2 n}, F_{ - 2 n} )$ directly inside $L(m, 0)
\otimes F_ { - 2 n ( n m +1)}$. This approach requires many deep
results on the structure of the vertex operator algebra $L(m,0)$
and deserves to be investigated independently.  Instead of this
approach, we realize the vertex algebra  $D_{m, 2 n} \otimes F_{ -
2 n} $ inside  a larger lattice vertex algebra. Then the formulas
for the generators are much simpler (cf. Section \ref{parne}). The
similar analysis can be done when $k$ is odd (cf. Section
\ref{neparne}). This approach was also used in \cite{A3} for
studying the fusion rules for the $N=2$ vertex operator
superalgebra $D_{m,1}$.

\section{Lattice and affine vertex algebras } \label{stwist}

In this section, we shall recall the lattice construction of
vertex superalgebras from \cite{DL}, \cite{K}.

 Let $L$  be a  lattice.
  Set ${\h}={\C}\otimes_{\Z}L$ and
extend the ${\Z}$-form $ \la \cdot, \cdot \ra $ on $L$ to ${\h}$.
 Let $\hat{{\h}}={\C}[t,t^{-1}]\otimes {\h} \oplus {\C}c$ be the affinization of
${\h}.$ We also use the notation $h(n)=t^{n}\otimes h$ for $h\in
{\h}, n\in {\Z}$.

Set
$
\hat{{\h}}^{+}=t{\C}[t]\otimes
{\h};\;\;\hat{{\h}}^{-}=t^{-1}{\C}[t^{-1}]\otimes {\h}.
$
Then $\hat{{\h}}^{+}$ and $\hat{{\h}}^{-}$ are abelian subalgebras
of $\hat{{\h}}$. Let $U(\hat{{\h}}^{-})=S(\hat{{\h}}^{-})$ be the
universal enveloping algebra of $\hat{{\h}}^{-}$. Let ${\l} \in
{\h}$. Consider the induced $\hat{{\h}}$-module
\begin{eqnarray*}
M(1,{\l})=U(\hat{{\h}})\otimes _{U({\C}[t]\otimes {\h}\oplus
{\C}c)}{\C}_{\l}\simeq
S(\hat{{\h}}^{-})\;\;\mbox{(linearly)},\end{eqnarray*} where
$t{\C}[t]\otimes {\h}$ acts trivially on ${\C}$,
${\h}$ acting as $\la h, {\l} \ra$ for $h \in {\h}$
and $c$ acts on ${\C}$ as multiplication by 1. We shall write
$M(1)$ for $M(1,0)$.
 For $h\in {\h}$ and $n \in {\Z}$ write $h(n) =  t^{n} \otimes h$. Set
$
h(z)=\sum _{n\in {\Z}}h(n)z^{-n-1}.
$

Then $M(1)$ is a vertex operator algebra which is generated by the
fields $h(z)$, $h \in {\h}$, and $M(1,{\l})$, for $\l \in {\h}$,
are irreducible modules for $M(1)$.

Let $\hat{L}$ be the canonical central extension of $L$ by the
cyclic group $\< \pm 1\>$:
\begin{eqnarray}\label{2.7}
1\;\rightarrow \< \pm 1\>\;\rightarrow \hat{L}\;\bar{\rightarrow}
L\;\rightarrow 1
\end{eqnarray}
with the commutator map $c(\alpha,\beta)=(-1)^{\< \alpha,\beta\>}$
for $\alpha,\beta \in L$. Let $e:  L\to \hat L$ be a section such
that $e_0=1$ and $\epsilon: L\times L\to \<\pm 1\>$ be the
corresponding 2-cocycle. Then
$\epsilon(\a,\b)\epsilon(\b,\a)=(-1)^{\<\a,\b\>},$
\begin{equation}\label{2c}
\e(\a,\b)\e(\a+\b,\gamma)=\e(\b,\gamma)\e(\a,\b+\gamma)
\end{equation}
and $e_{\a}e_{\b}= \e(\a,\b)e_{\a+\b}$ for $\a,\b,\gamma\in L.$
   Form the induced $\hat{L}$-module
\begin{eqnarray*}
{\C}\{L\}={\C}[\hat{L}]\otimes _{\< \pm 1\>}{\C}\simeq
{\C}[L]\;\;\mbox{(linearly)},\end{eqnarray*} where ${\C}[\cdot]$
denotes the group algebra and $-1$ acts on ${\C}$ as
multiplication by $-1$. For $a\in \hat{L}$, write $\iota (a)$ for
$a\otimes 1$ in ${\C}\{L\}$. Then the action of $\hat{L}$ on ${\C}
\{L\}$ is given by: $a\cdot \iota (b)=\iota (ab)$ and $(-1)\cdot
\iota (b)=-\iota (b)$ for $a,b\in \hat{L}$.

Furthermore we define an action of ${\h}$ on ${\C}\{L\}$ by:
$h\cdot \iota (a)=\< h,\bar{a}\> \iota (a)$ for $h\in {\h},a\in
\hat{L}$. Define $z^{h}\cdot \iota (a)=z^{\< h,\bar{a}\> }\iota
(a)$.

The untwisted space associated with $L$ is defined to be
\begin{eqnarray*}
V_{L}={\C}\{L\}\otimes _{{\C}}M(1)\simeq {\C}[L]\otimes S(\hat{\h}
^{-})\;\;\mbox{(linearly)}.
\end{eqnarray*} Then
$\hat{L},\hat{{\h}},z^{h}\;(h\in {\h})$ act naturally on $V_{L}$
by acting on either ${\C}\{L\}$ or $M(1)$ as indicated above.
Define ${\1}= \iota ( e_0) \in V_L$.
We use a normal ordering procedure, indicated by open colons,
which signify that in the enclosed expression, all creation
operators $h(n)$ $(n<0)$,$a\in \hat{L}$ are to be placed to the
left of all annihilation operators $h(n),z^{h}\;(h\in {\h},n\ge
0)$. For $a \in \hat{L}$, set
\begin{eqnarray*}
Y(\iota (a),z)= : e^{\int
(\bar{a}(z)-\bar{a}(0)z^{-1})}az^{\bar{a}}:.
\end{eqnarray*}
Let $a\in \hat{L};\;h_{1},\cdots,h_{k}\in
{\h};n_{1},\cdots,n_{k}\in {\Z}\;(n_{i}> 0)$. Set
\begin{eqnarray*}
v= \iota (a)\otimes h_{1}(-n_{1})\cdots h_{k}(-n_{k})\in
V_{L}.\end{eqnarray*}

 Define vertex operator $Y(v,z)$ with
\bea \label{defvertex}  :\left({1\over (n_{1}-1)!}({d\over
dz})^{n_{1}-1}h_{1}(z)\right)\cdots \left({1\over
(n_{k}-1)!}({d\over dz})^{n_{k}-1}h_{k}(z)\right)Y(\iota (a),z): .
\eea
 This gives us a well-defined linear map
\begin{eqnarray*}
Y(\cdot,z):& &V_{L}\rightarrow
(\mbox{End}V_{L})[[z,z^{-1}]]\nonumber\\ & &v\mapsto Y(v,z)=\sum
_{n\in {\Z}}v_{n}z^{-n-1},\;(v_{n}\in {\rm End}V_{L}).
\end{eqnarray*}

Let $\{\; h_{i}\;|\;i=1,\cdots,d\}$ be an orthonormal basis of
${\h}$ and set
\begin{eqnarray*}
\omega ={1\over 2}\sum _{i=1}^{d} h_{i}(-1) h_{i}(-1)\in V_{L}.
\end{eqnarray*}
Then $Y(\omega,z)=\sum_{n\in {\Z}}L(n) z^{-n-2}$ gives rise to a
representation of the Virasoro algebra on $V_{L}$ with the central
charged $d$ and
\begin{eqnarray}  \label{vir.rel}
& &L(0)\left(\iota(a)\otimes h_{1}(-n_{1})\cdots
h_{n}(-n_{k})\right)\nonumber \\ &=&\left({1\over 2}\<
\bar{a},\bar{a}\>+n_{1}+\cdots+n_{k}\right) \left(\iota(a)\otimes
h_{1}(-n_{1})\cdots h_{k}(-n_{k})\right).
\end{eqnarray}

The following theorem was proved in \cite{DL} and \cite{K}.
\begin{theorem}
The structure $(V_L, Y,{\1}, L(-1) )$ is a vertex superalgebra.
\end{theorem}

 Define the Schur polynomials $p_{r}(x_{1},x_{2},\cdots)$
 in variables $x_{1},x_{2},\cdots$ by the following equation:
\begin{eqnarray}\label{eschurd}
\exp \left(\sum_{n= 1}^{\infty}\frac{x_{n}}{n}y^{n}\right)
=\sum_{r=0}^{\infty}p_{r}(x_1,x_2,\cdots)y^{r}.
\end{eqnarray}
For any monomial $x_{1}^{n_{1}}x_{2}^{n_{2}}\cdots x_{r}^{n_{r}}$
we have an element $h(-1)^{n_{1}}h(-2)^{n_{2}}\cdots
h(-r)^{n_{r}}{\1}$ in both $M(1)$ and $V_{L}$ for $h\in{\h}.$
 Then for any polynomial $f(x_{1},x_{2}, \cdots)$,  $f(h(-1),
h(-2),\cdots){\1}$ is a well-defined element in $M(1)$ and
$V_{L}$.  In particular, $p_{r}(h(-1),h(-2),\cdots){\1}$  for $r
\in {\N}$  are elements of $M(1)$ and $V_{L}$.

Suppose $a,b\in \hat{L}$ such that $\bar{a}=\alpha,\bar{b}=\beta$.
Then
\begin{eqnarray}\label{erelation}
Y(\iota(a),z)\iota(b)&=&z^{\<\alpha,\beta\>}\exp\left(\sum_{n=1}^{\infty}
\frac{\alpha(-n)}{n}z^{n}\right)\iota(ab)\nonumber\\
&=&\sum_{r=0}^{\infty}p_{r}(\alpha(-1),\alpha(-2),\cdots)\iota(ab)
z^{r+\<\alpha,\beta\>}.
\end{eqnarray}
Thus
\begin{eqnarray}\label{4.7}
\iota(a)_{i}\iota(b)=0\;\;\;\mbox{ for }i\ge -\<\alpha,\beta\>.
\end{eqnarray}
Especially, if $\<\alpha,\beta\>\ge 0$, we have
$\iota(a)_{i}\iota(b)=0$ for  $i\ge 0$, and if
$\<\alpha,\beta\>=-n<0$, we get
\begin{eqnarray}\label{4.8}
\iota(a)_{i-1}\iota(b)=p_{n-i}(\alpha(-1),\alpha(-2),\cdots)\iota(ab)
\;\;\;\mbox{ for }i\in \{ 0, \dots, n\}.
\end{eqnarray}

 Let $n \in {\Z}$, and $\la \beta , \beta \ra= n$.
Define
$$L_n= {\Z} {\beta}, \quad F_n= V_{ L_n}. $$

Then $F_n$ is a vertex algebra if $n$ is even, and a vertex
superalgebra if $n$ is odd. For $i \in {\Z}$, let $\overline{i} =
i + n{\Z} \in  {\Z} / { n \Z}$. We define $ F_n ^{\overline{i}} =
V_{ \Z \beta + \frac{i}{n} \beta}$. Clearly $F_n = F_n
^{\overline{0}}$. It is well-known (cf. \cite{D}, \cite{DL},
\cite{Xu}) that the set $\{ F_n ^{\overline{i}} \} _{
i=0,\dots,\vert n \vert -1}$ provides all irreducible
$F_n$--modules. In particular, $F_n$ has $\vert n \vert$
inequivalent irreducible modules.

The fusion algebra is (cf. \cite{DL})
\bea \label{fa1} F_n ^{\overline{i}} \times F_n ^{\overline{j}} =
F_n ^{\overline{i+j} }. \eea
If $n=2k$ is even, we define
$ \tilde{L}_{2k}   = \frac{\beta}{2} + {\Z}{\beta},$ and $MF_{2k}
= V_{ \tilde{L}_{2k}} = F_{2k} ^{\overline{k}}$. Then $F_{2k}$ is
a vertex algebra, and $MF_{2k}$ is a $F_{2k} $--module.

We shall also need the following result from \cite{DLM}.

\begin{proposition} \label{regular2}
\cite{DLM} The vertex (super)algebra $F_n$ is  regular, i.e. any
(weak) $F_n$--module is completely reducible.
\end{proposition}
\vskip 10mm Let ${\g}$ be the Lie algebra $sl_2$ with generators
$e, f, h$ and relations $[e,f]=h$, $[h,e] = 2 e$, $[h, f ] =-2f$.
Let ${\hg} = {\g} \otimes {\C}[t, t ^{-1}] \oplus {\C} K$ be the
corresponding affine Lie algebra of type $A_1 ^{(1)}$. As usual we
write $x(n)$ for $x \otimes t ^{n}$ where $x \in {\g}$ and $n \in
{\Z}$. Let $\Lambda_0$, $\Lambda_1$ denote the fundamental weights
for $\hg$. For any complex numbers $m,j$, let $L(m,j) = L( (m-j)
\Lambda_0 + j \Lambda_1)$ be the irreducible highest weight
$\hat{sl_2}$--modules with the highest weight $(m-j) \Lambda_0 + j
\Lambda_1$ . Then $L(m,0)$ has a natural structure of a simple
vertex operator algebra. Let  ${\1} _m$ denote the vacuum vector
in $L(m,0)$.

 If $m$ is a positive integer then $L(m,0)$ is a regular vertex operator algebra,
  and the set $\{ L(m,j)
\}_{j=0,\dots,m }$ provides all inequivalent irreducible
$L(m,0)$--modules. The fusion algebra (cf. \cite{FZ}) is given by
\bea \label{fa2} L(m,j) \times L(m,k) = \sum_{ i=\mbox{max} \{0,
j+k -m \} } ^{\mbox{min} \{j, k\} } L(m, j+ k - 2i). \eea
In particular, $L(m,m) \times L(m,j) = L(m,m-j)$.

We shall now recall the lattice construction of the vertex
operator algebra $L(m,0)$. Define the following lattice
\bea && A_{1,m}= {\Z} {\alpha}_1 + \cdots + {\Z} {\alpha}_m
\nonumber \\
 && \langle\alpha_i, \alpha_j \rangle =2 \delta_{i,j},  \nonumber \eea
 for every $i, j \in \{ 1, \dots, m \}. $
Define also $ \tilde{A}_{1,m} = \frac{\alpha _1 + \cdots +
\alpha_m }{2} + A_{1,m}$. We have:

\begin{lemma} \cite{DL} \label{2.1}
The vectors
${E} = \iota( e_{ \alpha_1} ) + \cdots +  \iota( e_{ \alpha_m} )$,
${F} = \iota( e_{ -\alpha_1} ) + \cdots + \iota( e_{- \alpha_m}
)$, span a subalgebra of $V_{ A_{1,m}}$ isomorphic to $L(m,0)$.
Moreover, $L(m,m)$ is a $L(m,0)$ submodule of $V_{ \tilde{A}_{1,m}
}$.
\end{lemma}

\section{ The definition of $D_{m,k}$ } \label{def-d}

In this section we give the definition of the vertex operator
(super)algebra $D_{m ,k}$. Let the vertex (super)algebras $L(m,0)$
and $F_k$ be defined as in Section \ref{stwist}.

\begin{definition}
Let $m \in {\C} \setminus \{ 0, -2 \}$, and let $k$ be a
nonnegative integer. Then $D_{m,k}$ is a vertex subalgebra of the
vertex operator (super)algebra $L(m,0) \otimes F_{k}$ generated by
the vectors :
$$
X= e(-1) {\1}_m \otimes \iota(e_{ \beta}) \quad \mbox{ and} \quad
Y= f(-1) {\1}_m \otimes \iota(e_{- \beta}). $$
\end{definition}
\vskip 5mm  Let ${\1}_{m,k} = {\1}_m \otimes {\1} \in D_{m,k}
\subset L(m,0) \otimes F_k$. Define also
\bea
&& H= h(-1) {\1}_m \otimes {\1} + {\1} \otimes \beta(-1) {\1},
\nonumber \\
&&
 \omega_{ m , k} = \frac{1}{ 2 (m +
2)} \left(
X_{ k - 1} Y + Y_{ k -1} X + \frac{1 - k}{   m k + 2} H_{-1} ^{2}
{\1}_{m,k} \right) . \nonumber \eea
It is easy to see that the components of the field
$$Y(\omega_{m,k},z) = \sum_{ n \in {\Z} } L(n) z ^{- n -2}$$
give rise a representation of the Virasoro algebra of central
charge $c_m = \frac{3 m}{m +2}$.  For $n \ge 0$ one has
$$ L(n) x = \delta_{n,0} ( 1 + \frac{k}{2}) x \quad  \mbox{and}
\quad L(n) y = \delta_{n,0} ( 1 + \frac{k}{2}) y.$$
Thus, $L(0)$ defines on $D_{m,k}$ a $ {\Z}_+$--graduation if $k$
is even, and a $\frac{1}{2} {\Z}_+$--graduation if $k$ is odd. In
this way we get the following theorem.

\begin{theorem} Let $m \in {\C} \setminus \{ 0, -2 \}$,  and let $k$ be a
nonnegative integer. Then $D_{m,k}$ is a vertex operator algebra
if $k$ is even and a vertex operator superalgebra if $k$ is odd.
The Virasoro vector is $\omega_{m,k}$, the vacuum vector is
${\1}_{m,k}$ and the rank is $c_m$.
\end{theorem}

\begin{remark} Let $k=0$. Then $D_{ m, 0 }$ is ismorphic to the
$\hat {sl_2}$ vertex operator algebra $L(m,0)$.
 Note also that the vector
$$\omega_{ m , 0} = \frac{1}{ 2 (m + 2)} \left(
X_{ - 1} Y + Y_{ -1} X + \frac{1 }{ 2 } H_{-1} ^{2} {\1}_{m,0}
\right)$$
 coincides with the Virasoro vector in $L(m,0)$
constructed using the Sugawara construction.
\end{remark}

\begin{remark}
 For $k=1$,
$D_{m,1}$ is in fact the vertex operator superalgebra associated
to the vacuum representation  of the $N=2$ superconformal algebra
constructed using the Kazama-Suzuki mapping (cf. \cite{FST1},
\cite{KS}). The Virasoro vector in $D_{m,1}$ is
$$\omega_{ m , 1} =  \frac{1}{ 2 (m + 2)}(
X_{ 0} Y + Y_{ 0} X  ).$$
Its representation theory  was studied in \cite{EG} and \cite{A2}.
In particular, when $m$ is not a nonnegative integer then
$D_{m,1}$ is not rational. In Theorem \ref{generalizacija} we will
generalize this fact for every positive integer $k$.

It was proved in  \cite{A3} that if $m$ is a positive integer,
then $D_{m,1}$ is a regular vertex operator superalgebra and that
the vertex superalgebra $D_{m,1} \otimes F_{-1}$ is a simple
current extension of the vertex algebra $L(m,0) \otimes F_{ - 2
(m+2)}$.
\end{remark}

The definition of $D_{m,k}$ implies that for every weak
$L(m,0)$--module $M$, $M \otimes F_{k}$ is a weak module for
$D_{m,k}$. Thus, the representation theory of $D_{m,k}$ is closely
related to the representation theory of the vertex operator
algebra $L(m,0)$. The case when $m$ is an nonnegative integer will
be studied in following sections. When $m \ne -2$ and $m$ is not
an admissible rational number, then every highest weight $
\hat{sl_2}$--module of level $m$ is a module for the vertex
operator algebra $L(m,0)$. This easily gives that $D_{m,k}$ is not
rational.
 In the case when $m$ is an
admissible rational number,
 using the similar arguments as in \cite{A2}, and
using the representation theory of the vertex operator algebra
$L(m,0)$ in this case (cf. \cite{AM}) one can construct infinitely
many inequivalent irreducible $D_{m,k}$--modules. In order to be
more precise, we shall state the following lemma.


\begin{lemma} \label{konst-mod} Assume that
$m$ is not a nonnegative integer and $m \ne -2$. Let $k \ge 1$.
 Then for every $t \in {\C}$ there is an ordinary $D_{m,k}$--module $N_t$
 such that $N_{t} = \oplus_{ n \in \frac{1}{2} {\Z}_+ }
N_{t} (n)$, and  the top level $N_t (0)$ satisfies $$  N_{t} (0) =
{\C} w, \quad L(n) w = t \delta_{n,0} w \ \ \mbox{for} \ n \ge 0.
$$
\end{lemma}
{\em Proof.} The proof will use a similar consideration to those
in \cite{A2}, Section 6.

Assume that $ m$ is not a positive integer and $t\in {\C}$. The
results from \cite{AM} give that for every $ q \in {\C}$ there is
a ${\Z}_+$--graded $L(m,0)$--module $M_{q} = \oplus_{ n \in
{\Z}_+} M_{q} (n)$ and a weight vector $v_q \in M_{q} (0)$ such
that
\bea && \Omega (0) \vert M_{q} (0)   \equiv \frac{(m+2) m}{2}
\mbox{Id}, \quad
h(0) v_{q} = q v_{q}, \nonumber \eea
where $\Omega (0) =  e(0) f(0) + f(0) e(0) + \frac{1}{2} h(0)
^{2}$ is the Casimir element acting on the $sl_2$--module $M_{q}
(0)$.
Then $M_{q} \otimes F_k$ is a weak $D_{m,k}$--module. Choose $q
\in {\C}$ such that
$$ \frac{m}{4} - \frac{1}{4 (m+2)} q ^{2} + \frac{1-k}{2( m+2) (m
k + 2) } q ^{2} = t. $$
Let $N_t$ be the $D_{m,k}$--submodule of $M_q \otimes F_k$
generated by the vector $w = v_q \otimes {\1}$. Then for $n \ge 0$
we have that
$$ L(n) w = \delta_{n,0} (\frac{m}{4} - \frac{1}{4 (m+2)} q ^{2} +
\frac{1-k}{2( m+2) (m k + 2) } q ^{2} ) w = \delta_{n,0} t  w. $$
Now it is easy to see  that $N_t$ is an ordinary
$\frac{1}{2}{\Z}_+$--graded $D_{m,k}$--module with the top level $
N_t (0) = {\C} w$ and that $ L (0) \vert N_t (0) \equiv t
\mbox{Id}$. Thus, the lemma holds. \qed

In fact, Lemma \ref{konst-mod} gives that there is uncountably
many inequivalent irreducible $D_{m,k}$--modules. Thus, we
conclude that   the following theorem holds.

\begin{theorem} \label{generalizacija}
Assume that $m$ is not a nonnegative integer and $m \ne -2$. Then
for every positive integer $k$, the vertex operator (super)algebra
$D_{m,k}$ is not rational.
\end{theorem}

\begin{remark} In what follows we will prove that if $m$ is a positive integer,
then $D_{m,k}$ is rational. In fact, we will establish more
general complete reducibility theorem, which will imply that
$D_{m,k}$ is regular in the sense of \cite{DLM}.
\end{remark}

\section{ The lattice construction of  $D_{ m, k}$ for $ m \in {\N}$ }
\label{construction}

 In this section we give the lattice construction of the vertex
 operator (super)algebra $D_{m,k}$. This construction is a
 generalization of the lattice constructions of the vertex operator
 algebra $L(m,0)$ (cf. \cite{DL} and our Lemma \ref{2.1}) and of
 the N=2 vertex operator superalgebra $L_{ c_m}$ (cf. \cite{A3}).

Let $m ,k$ be positive integers. Define the lattice
 \bea && \Gamma_{m,k}= {\Z} {\gamma}_1 + \cdots + {\Z} {\gamma}_m \nonumber \\
 && \langle \gamma_i, \gamma_j \rangle =2 \delta_{i,j} + k,   \nonumber \eea
 for every $i, j \in \{ 1, \dots, m \}. $

 Then $V_{ \Gamma_{m,k}}$ is a vertex operator algebra if $k$ is even
 and a vertex operator superalgebra if $k$ is add.

\begin{proposition} Let $m , k $ be positive integers. The vertex operator (super)algebra
 $D_{m,  k}$ is isomorphic to the subalgebra of the vertex operator
(super)algebra $V_{\Gamma_{m,k}}$ generated by the vectors
\bea && \bar{X} = \iota(e_{\gamma_1}) + \cdots
+\iota(e_{\gamma_m}); \nonumber
\\
&& \bar{Y} = \iota(e_{-\gamma_1}) + \cdots +\iota(e_{-\gamma_m}).
\nonumber
  \eea
   Set $\bar{H} = \bar{X}_k \bar{Y}$. Then the  Virasoro vector in
$D_{m,k}$ is given by

\bea  \bar{\omega}_{ m , k} &&=  \frac{1}{ 2 (m + 2)} \left(
{\bar X}_{ k - 1} {\bar Y} + {\bar Y}_{ k -1} {\bar X } + \frac{1
- k}{ m k + 2} {\bar H}_{-1} ^{2} {\1} \right)  \nonumber \\&&=
\frac{1}{2 ( m + 2)} \sum_{ i = 1} ^{m} \gamma_{ i} (-1 ) ^{2}
{\1} + \frac{1}{ m+2} \sum_{ i \ne j} \iota ( e_{ \gamma_i -
\gamma _j} ) + \nonumber \\
&&+ \frac{1 - k}{ 2 ( m + 2) ( m k + 2) } \left(\sum_{ i = 1 }
^{m} \gamma_i (-1) \right) ^{2} {\1}. \nonumber \eea
\end{proposition}
{\em Proof.} Define the lattice $\Gamma_1$ by
 \bea && \Gamma_1 = {\Z} {\alpha_1} + \cdots + {\Z} {\alpha_m} +
 {\Z} \beta, \nonumber \\
 && \langle \alpha_i, \alpha_j\rangle = 2 \delta_{i,j}, \ \
  \langle \alpha_i, \beta \rangle =
 0, \ \ \langle \beta, \beta \rangle = k. \nonumber \eea
For $i=1, \dots, m$ set $\gamma_i = \alpha_i + \beta$.
It is clear that the lattice $\Gamma_{m,k}$ can be identified with
the sublattice ${\Z}{\gamma_1} + \cdots + {\Z} {\gamma_m}$ of the
lattice $\Gamma_1$. In the same way $V_{ \Gamma_{m,k}}$ can
treated as a subalgebra of the vertex operator (super)algebra $V_{
\Gamma_1}$.
Lemma \ref{2.1} implies that
${\bar E} = \iota( e_{ \alpha_1} ) + \cdots +  \iota( e_{
\alpha_m} )$, ${\bar F} = \iota( e_{ -\alpha_1} ) + \cdots +
\iota( e_{- \alpha_m} )$, span a subalgebra of $V_{ \Gamma_1}$
isomorphic to $L(m,0)$, and the elements $\iota( e_{ \beta})$,
$\iota( e_{ -\beta})$ span a subalgebra isomorphic to $F_k$. Since
$$ \bar{X} = {\bar E}_{ -1} \iota( e_{ \beta}) \ \ \mbox{and} \ \
\bar{Y} = {\bar F}_{ -1}
%
 \iota( e_{ -\beta}), $$
we conclude that the vertex subalgebra generated by the elements
$\bar{X}, \bar{Y}  \in V_{ \Gamma_{m,k}}  \subset V_{ \Gamma_1}$
is isomorphic to the vertex operator (super)algebra $D_{m,k}$.
This concludes the proof of the theorem. \qed

The previous result implies that we can identify the generators of
$D_{m,k}$ in $L(m,0) \otimes F_{k}$ with the generators of
$D_{m,k}$ in $V_{ \Gamma_{m,k}}$. Thus we can assume that $X
\approx\bar{X}$, $Y \approx \bar{Y}$. This identification implies
that $H \approx \bar{H}$ and
$\omega_{m,k}\approx\bar{\omega}_{m,k}$.

We shall also prove an interesting proposition which identifies
some regular subalgebras of $D_{m,k}$.

\begin{proposition} \label{reg-sub}  For every positive integer $n$ we have that
\bea
&\iota( e_{   n ( \gamma_1 + \dots +\gamma_m ) } ), \ \iota( e_{ -
n  ( \gamma_1 + \dots +\gamma_m ) } ) \in D_{ m, k}. \nonumber
\eea
In particular, $D_{m , k}$ has a vertex subalgebra isomorphic to
$F_{ n ^{2} m ( m k + 2)  }$.
\end{proposition}

{\em Proof.} Using relations (\ref{4.7}) and (\ref{4.8}), it is
easy to prove  that:
\bea
&& \left( X_{ -( n m - 1) k - 2 n + 1}  \cdots X_{ -  (n -1) m - 2
n + 1} \right)
 \cdots
 \left(
 X_{ -(2 m -1) k - 3} \cdots  X_{ - m k -3}
 \right)
 \cdot \nonumber
 \\ && \cdot
  \left(
  X_{ -(m-1) k - 1} \cdots X_{ -k - 1} X_{ -1}
 {\1}
 \right)
  =C \iota( e_{ n ( \gamma_1 + \dots +\gamma_m ) } ) \nonumber
\eea
for some constant $C$. Thus $\iota( e_{   n ( \gamma_1 + \dots
+\gamma_m ) } ) \in D_{m,k}$. Similarly we prove that \\ $\iota(
e_{ -n ( \gamma_1 + \dots +\gamma_m ) } ) \in D_{m,k}$. The second
assertion of the proposition follows from the fact that the
vectors
 $\iota(
e_{ \pm n ( \gamma_1 + \dots +\gamma_m ) } )$ generate a
subalgebra of $V_{\Gamma_{m,k}} $ isomorphic to $F_{ n ^{2} m (m k
+2)}$.
 \qed

 \section{  Regularity of the vertex operator algebra $D_{m, 2 n}$}
\label{parne}

 In this section we study the vertex algebra $L(m,0) \otimes F_{
 -2 n (m n +1)}$ where $m,n$ are positive integers.  We now that $L(m,0) \otimes F_{
 -2 n (m n+1)}$ is a regular vertex algebra; i.e. every module for this
 vertex algebra is completely reducible. Its irreducible modules
 are :

 $$ L(m,r) \otimes F_{ -2 n ( m n+1)} ^{\bar{s}} , \quad r \in \{1, \dots,
 m\}, \ \bar{s} \in \frac{ \Z} { -2 n (m n + 1) {\Z} }. $$

 The fusion rules can be calculated easily from the fusion rules
 for $L(m,0)$ and $F_{-2 n (m n +1)}$.

 Our main goal is to show that the vertex operator algebra $D_{m, 2 n}$ is
 isomorphic to a subalgebra  of $L(m,0)
 \otimes F_{-2 n (m n +1)}$. In order to do this, we shall  first
 give the lattice construction of the vertex algebra $L(m,0)
 \otimes F_{-2 n (m n +1)}$.

 Define the following lattice :
 \bea && L= {\Z} {\alpha}_1 + \cdots + {\Z} {\alpha}_m + {\Z}
 \beta  \nonumber \\
 && \langle \alpha_i, \alpha_j \rangle =2 \delta_{i,j}, \quad
 \langle \alpha_i, \beta
 \rangle =0, \quad \langle \beta, \beta \rangle = - 2 n ( m n+1) \nonumber \eea
 for every $i, j \in \{ 1, \dots, m \}. $

We shall now give another description of the lattice $L$.

 For $i=1, \dots, m$, we define
 \bea && \delta =  n {\alpha}_1 + \cdots + n {\alpha_m} + \beta ,
 \nonumber \\
 && {\gamma}_i = {\alpha}_i + \delta. \nonumber \eea
 Since
 $$ \alpha_i =  \gamma_i - \delta, \  \beta = ( n m+1) \delta -
n (\gamma_1 + \cdots + \gamma_m),$$
we have that
  \bea && L= {\Z} {\gamma}_1 + \cdots + {\Z} {\gamma}_m + {\Z}
 \delta  \nonumber \\
 && \langle \gamma_i, \gamma_j \rangle =2 \delta_{i,j} + 2 n , \quad
 \langle \gamma_i, \delta
 \rangle
 =0, \quad \langle \delta, \delta \rangle = - 2 n  \nonumber \eea
 for every $i, j \in \{ 1, \dots, m \}. $

 In fact, we have proved that
\bea \label{parne-rel}
 && L \cong  \Gamma_{m,2 n} + L_{- 2 n} \cong  A_{1,m} + L_{ -
2 n ( m n + 1) }, \eea
which implies that \bea \label{parne-rel2}
 && V_L \cong  V_{ \Gamma_{m,2 n}}  \otimes F_{- 2 n} \cong  V_{ A_{1,m}}
  \otimes F_{ -
2 n ( m n + 1) }. \eea

Define the following vectors in the vertex algebra $V_L$:
\bea && E = \iota(e_{\alpha_1}) + \cdots +\iota(e_{\alpha_m});
\nonumber \\
&& F = \iota(e_{-\alpha_1}) + \cdots +\iota(e_{-\alpha_m}).
\nonumber  \eea
These vectors span a subalgebra of $V_L$ isomorphic to $L(m,0)$.

As in Section \ref{construction} we define:
\bea && X = \iota(e_{\gamma_1}) + \cdots +\iota(e_{\gamma_m});
\nonumber
\\
&& Y = \iota(e_{-\gamma_1}) + \cdots +\iota(e_{-\gamma_m}).
\nonumber
  \eea
Clearly $X,Y$ span a subalgebra isomorphic to $D_{m, 2 n}$. In
fact, the definition of elements $E$, $F$, $X$, $Y$ together with
relations (\ref{parne-rel}) and (\ref{parne-rel2}) imply the
following lemma.

\begin{lemma} \label{sub1} \item[(1)] Let $V$ be the subalgebra of $V_L$ generated by
the vectors $$E, F, \iota(e_{\beta}), \iota(e_{-\beta}). $$ Then
 $V \cong L(m,0) \otimes F_{-2 n (m n +1)}$.

\item[(2)] Let $W$ be the subalgebra of $V_L$ generated by the vectors
 $$X, Y , \iota(e_{\delta}),\iota(e_{-\delta}). $$
  Then  $ W\cong D_{m, 2 n} \otimes F_{-2 n}.$ \end{lemma} \vskip 5mm

Now using standard calculations in lattice vertex algebras one
easily gets the following  important lemma.

\begin{lemma} \label{relations} In the vertex algebra $V_L$ the following relations
hold:
\item[(1)] $X = (E_{-2 n -1}  \iota( e_{n( \alpha_1 + \cdots +\alpha_m)})
)_{-1} \iota ( e_{\beta})$;

\item[(2)] $Y = (F_{-2 n -1}  \iota( e_{ - n (\alpha_1 + \cdots +\alpha_m)})
) _{-1} \iota(e_{-{\beta} })$;

\item[(3)] $ \iota (e _{\delta}) = \iota( e_{ n (\alpha_1 + \cdots
+\alpha_m)} ) _{-1} \iota(e_{\beta})$;

\item[(4)] $ \iota (e _{-\delta}) = \iota( e_{ - n (\alpha_1 + \cdots
+\alpha_m)} ) _{-1} \iota(e_{-\beta})$;

\item[(5)] $E = X_{-1} \iota (e_{-\delta}) $  ;

\item[(6)] $F = Y_{-1} \iota (e_{\delta}) $  ;

\item[(7)] $   \iota (e_{\beta}) =  \iota (e_{ (n m+1)\delta}) _{-1}
  \iota (e_{-  n (\gamma_1 + \cdots + \gamma_m )})  $  ;

\item[(8)] $   \iota (e_{-\beta}) =  \iota (e_{ -(n m+1)\delta}) _{-1}
  \iota (e_{n({\gamma_1 + \cdots +\gamma_m })})  $  .
\end{lemma}

\vskip 5mm

\begin{theorem} \label{isom1} The vertex subalgebras $V$ and $W$ coincide. In
particular, we have the following isomorphism of vertex algebras:
\bea L(m,0) \otimes F_{ -2 n (m n +1)} \cong D_{m, 2 n} \otimes
F_{-2 n } . \nonumber \eea
\end{theorem}
{\em Proof.} Using the same arguments as in the proof of
Proposition \ref{reg-sub} we get
$$ \iota(e_{ \pm n ( \alpha_1 + \cdots + \alpha_m ) } ) \in V,
\quad \iota(e_{ \pm n ( \gamma_1 + \cdots + \gamma_m) } ) \in W.
$$

Then the relations (1) - (4) in Lemma \ref{relations} implies that
$X, Y, \iota(e_{ \pm \delta}) \in V$. Thus $W \subset V$.
Similarly, the relations (5) - (8) in Lemma \ref{relations} gives
that $V \subset W$. Hence, $V=W$. Then Lemma \ref{sub1}   implies
that $ L(m,0) \otimes F_{ -2 n (m n+1)} \cong D_{m, 2 n}  \otimes
F_{-2 n} $. \qed

\vskip 5mm

The following proposition was essentially proved in \cite{DLM} and
\cite{DMZ}.

\begin{proposition}  \label{reg-dmz} Let $V$ be  a  vertex operator
(super) algebra and $s \in {\Z}$. Then $V \otimes F_s$ is a
regular vertex superalgebra if and only if $V$ is a regular vertex
operator (super)algebra.
\end{proposition}

\begin{theorem} Let $m, m_1,   \dots, m_r$ be positive integers
and let $k, k_1, \dots ,k_r$ be positive even integers.
\item[(1)] The vertex operator algebra $D_{m, k}$ is regular; i.e. every weak $D_{m, k}$--module
is completely reducible.
\item[(2)] The vertex operator algebra $D_{m_1, k_1} \otimes \cdots \otimes D_{ m_k,  k_r}$ is
regular.
\end{theorem}
{ \em Proof.} Since $L(m,0)$ and $F_{ - 2 n (n m +1)}$ are regular
vertex algebras, Proposition \ref{reg-dmz} implies that $L(m,0)
\otimes F_{ - 2 n ( n m +1)}$ is regular. Since
$$L(m,0) \otimes F_{ - 2 n ( n m +1)} \cong D_{m, 2 n} \otimes
F_{- 2n},$$ %
 using again Proposition \ref{reg-dmz} we get that $D_{m, 2 n}$ is a
regular vertex operator algebra. This gives (1). The proof of (2)
is now standard (cf. \cite{DLM}). \qed

Since $L(m,0)$ has $(m+1)$ inequivalent irreducible modules, and
for every $k \in {\Z}$ $F_k$ has $\vert k \vert$ inequivalent
irreducible modules , we get:
\begin{corollary} The vertex operator algebra $D_{m, 2 n}$ has exactly
$ (m+1) ( n m + 1)$ inequivalent  irreducible representations.
\end{corollary}

\section{ Regularity of the vertex operator superalgebra
$D_{m,k}$ for $k$ odd } \label{neparne}

In this section, we shall consider the case when $k$ is an odd
natural number. When $k=1$, then $D_{ m ,1}$ is the vertex
operator superalgebra associated to the
unitary vacuum representation for the $N=2$ superconformal
algebra. This case was studied in \cite{A3}.

First we see that the following relation between lattices holds:

\bea \label{nep-rel}
&& \Gamma_{m,k} + L_{-k} \cong ( A_{1,m} + L_{ - 2  k ( m k + 2)})
\cup  ( \tilde{A}_{1,m} + \tilde{L}_{ - 2  k ( m k + 2)}),
\eea
which implies the following isomorphism of vertex algebras:
\bea \label{nep-rel2}
&& V_{ \Gamma_{m,k} } \otimes F_{-k} \cong ( V_{ A_{1,m}}  \otimes
F_{ - 2 k ( m k + 2)}) \oplus  ( V_{ \tilde{A}_{1,m} } \otimes
MF_{ - 2 k ( m k + 2)}).
\eea

Using (\ref{nep-rel}), (\ref{nep-rel2})  and a completely
analogous proof to that of Theorem 7.1   in \cite{A3}, we get the
following result.

\begin{theorem}  \label{nep-main} We have the following isomorphism of  vertex
superalgebras:
$$ D_{m,k} \otimes F_{ -k} \cong L(m,0)  \otimes F_{ - 2k ( k m +2
)} \oplus L(m,m) \otimes MF_{ - 2 k (k m +2) }. $$
In other words, the vertex superalgebra $D_{m,k} \otimes F_{ -k}$
is a simple current extension of the vertex algebra  $L(m,0)
\otimes F_{ -2k ( k m +2 )}$.
\end{theorem}

Using Proposition \ref{reg-dmz}, Theorem \ref{nep-main} and the
fact that a simple current extension of a regular vertex algebra
is a regular vertex (super)algebra (cf. \cite{Li2}) we get the
following theorem.

\begin{theorem} Let $m, m_1,   \dots, m_r$ be positive integers
and let $k, k_1, \dots ,k_r$ be positive odd integers.
\item[(1)] The vertex operator superalgebra $D_{m, k}$ is regular.

\item[(2)] The vertex operator superalgebra $D_{m_1, k_1} \otimes \cdots \otimes D_{ m_k,  k_r}$ is
regular.
\end{theorem}

We also have:
\begin{corollary} The vertex operator superalgebra $D_{m,k}$ has
exactly $  \frac{ (m+1) ( k m +2)}{2}$ inequivalent irreducible
representations.
\end{corollary}
{\em Proof.} The results from \cite{Li2} imply that the extended
vertex superalgebra
 $$ L(m,0)  \otimes F_{ - 2k ( k m +2 )} \oplus
L(m,m) \otimes MF_{ - 2 k (k m +2) } $$
has exactly $ \frac{1}{2} (m+1) k (k m +2)$ inequivalent
irreducible representations (see also  \cite{A3}, \cite{Li3}).
Since the vertex superalgebra $F_{ -n}$ has $n$ inequivalent
irreducible representations, we conclude that $D_{m,k}$ has to
have $ \frac{ (m+1)  (k m +2)}{2}$ inequivalent irreducible
representations. \qed

\end{document}